%% file: agt-3-35.tex
\documentclass{gtart}
\input agtout

\lognumber{35}
\volumenumber{3}
\volumeyear{2003}
\papernumber{35}
\published{12 October 2003}
\pagenumbers{1037}{1050}
\received{17 Febuary 2003}
\revised{7 September 2003}
\accepted{10 October 2003}

\usepackage{amsmath,amssymb} 
\usepackage{enumerate}
\input xypic
\xyoption{all}

\theoremstyle{plain}
\newtheorem{theorem}{Theorem}[section]
\newtheorem{lemma}[theorem]{Lemma}

\newtheorem{proposition}[theorem]{Proposition}

\theoremstyle{definition}

\theoremstyle{remark}
\newtheorem{remark}[theorem]{Remark}

\makeatletter\let\c@equation=\c@theorem\makeatother

\newcommand{\Map}{\operatorname{Map}}
\newcommand{\Or}{{\rm Or}}

\newcommand{\res}{{\rm res}} 

\newcommand{\id}{\operatorname{id}}

\newcommand{\x}{\times}

\newcommand{\einsu}{[1,\infty)}

\newcommand{\point}{\mathit{pt}}

\newcommand{\cala}{{\mathcal{A}}}

\newcommand{\calc}{{\mathcal{C}}}
\newcommand{\cald}{{\mathcal{D}}}
\newcommand{\cale}{{\mathcal{E}}}
\newcommand{\calf}{{\mathcal{F}}}
\newcommand{\calg}{{\mathcal{G}}}

\newcommand{\cali}{{\mathcal{I}}}

\newcommand{\caln}{{\mathcal{N}}}

\newcommand{\calq}{{\mathcal{Q}}}

\newcommand{\cals}{{\mathcal{S}}}

\newcommand{\calv}{{\mathcal{V}}}

\newcommand{\calvc}{{\mathit{\calv\calc}}}
\newcommand{\calfin}{{\mathit{\calf\cali\caln}}}

\newcommand{\IA}{{\mathbb A}} 
\newcommand{\IB}{{\mathbb B}}

\newcommand{\IE}{{\mathbb E}} 
\newcommand{\IF}{{\mathbb F}}
 
\newcommand{\IH}{{\mathbb H}}

\newcommand{\IK}{{\mathbb K}}
\newcommand{\IL}{{\mathbb L}}

\newcommand{\IN}{{\mathbb N}}

\newcommand{\IR}{{\mathbb R}}

\newcommand{\IZ}{{\mathbb Z}}

\newcommand{\IKi}{{\IK^{-\infty}}}
\newcommand{\ILi}{{\IL^{-\infty}}}

\newcommand{\BFA}{{\bf A}}
\newcommand{\BFB}{{\bf B}}

\newcommand{\BFE}{{\bf E}}
\newcommand{\BFF}{{\bf F}}
\newcommand{\BFG}{{\bf G}}

\newcommand{\BFK}{{\bf K}}

\begin{document}
\title{On the domain of the assembly map\\in algebraic $K$--theory}
\asciititle{On the domain of the assembly map in algebraic K-theory}
\author{Arthur C. Bartels}
\address{SFB 478, Westf{\"a}lische Wilhelms-Universit{\"a}t\\48149 
M{\"u}nster, Germany}
\asciiaddress{SFB 478, Westfalische Wilhelms-Universitat\\48149 Munster, 
Germany}
\email{bartelsa@math.uni-muenster.de}


\begin{abstract}
We compare the domain of the assembly map in algebraic $K$--theory
with respect to the family of finite subgroups with the domain of 
the assembly map with respect to the family of virtually cyclic
subgroups and prove that the former is a direct summand of the later.
\end{abstract}
\asciiabstract{We compare the domain of the assembly map in algebraic 
K-theory with respect to the family of finite subgroups with the
domain of the assembly map with respect to the family of virtually
cyclic subgroups and prove that the former is a direct summand of the
later.}

\primaryclass{19D50} 

\secondaryclass{19A31, 19B28}

\keywords{$K$--theory, group rings, isomorphism conjecture}
\asciikeywords{K-theory, group rings, isomorphism conjecture}

\maketitle


\section{Introduction}

In algebraic $K$--theory assembly maps relate the algebraic $K$--theory of a 
group ring $R\Gamma$ to the algebraic $K$--theory of $R$ and the group 
homology of $\Gamma$. In the formulation of Davis and L{\"u}ck
\cite{Davis-Lueck(1998)} there is 
for every family of subgroups $\calf$ of $\Gamma$ an assembly map
\begin{eqnarray} \label{assemblymap}
H_*^{\Or \Gamma}(E\Gamma(\calf);\BFK R^{-\infty}) \to K_*(R\Gamma)
\end{eqnarray}
and these maps are natural with respect to inclusions of families of subgroups.
The notation is reviewed in more detail in Section~\ref{sec:orbit}.
The Isomorphism Conjecture of Farrell--Jones \cite{Farrell-Jones(1993)} for
algebraic $K$--theory (and $R = \IZ$) states that (\ref{assemblymap}) is an 
isomorphism, provided that $\calf = \calvc$ 
is the family of virtually cyclic subgroups.
This conjecture has been proven for  different classes of groups, 
cf. \cite{Farrell-Jones(1993)} \cite{Farrell-Jones(1998)}. 
Arbitrary coefficient rings are considered in 
\cite{Bartels-Farrell-Jones-Reich(2001)}. The assembly map is also
studied with $\calf = \calfin$ the family of finite subgroups or $\calf$ 
the family consisting of the trivial subgroup. 
For the trivial family there are
injectivity results for different classes of groups, 
cf.~\cite{Boekstedt-Hsiang-Madsen(1993)}, 
\cite{Carlsson-Pedersen(1995)}. Both results have been
extended to injectivity results for $\calf = \calfin$, see
\cite{Rosenthal(2002)} and
recent work of L{\"u}ck--Reich--Rognes--Varisco. 

In this paper we study the map 
\begin{eqnarray} \label{ktheorysplit}
H_*^{\Or \Gamma}(E\Gamma(\calfin);\BFK R^{-\infty})
\to
H_*^{\Or \Gamma}(E\Gamma(\calvc);\BFK R^{-\infty}).
\end{eqnarray}
It has been conjectured in \cite[p.260]{Farrell-Jones(1993)} (for $R=\IZ$) 
that this map
is  split injective. In various cases this follows from the above mentioned 
results. The purpose of this paper is to verify this conjecture in general.

\begin{theorem}   \label{main}
The map (\ref{ktheorysplit}) is split injective for 
arbitrary groups and rings.
\end{theorem}

In general the left hand side of (\ref{ktheorysplit}) is much better 
understood than the right hand side, cf.~\cite{Lueck(2000)}. 
Thus modulo the isomorphism conjecture Theorem~\ref{main} may be viewed as 
splitting a well understood factor from the $K$--theory of the group ring. 

For virtually cyclic groups Theorem~\ref{main} asserts that the
assembly map for the family $\calfin$ is split injective. This is 
a special case  of \cite{Rosenthal(2002)}. The language of 
$\Or \Gamma$--spectra from \cite{Davis-Lueck(1998)} allows us to extend
this splitting to the more general setting in (\ref{ktheorysplit}).

There is  a corresponding splitting result for $L$--theory: 
If we use $L^{-\infty}$--theory and $R$ and $\Gamma$ are such that 
$K_{-i}(RV) = 0$ for all virtually cyclic subgroups $V$ of $\Gamma$ and 
sufficiently large $i$, then (\ref{ktheorysplit})
remains split injective. This assumption is  satisfied if $R = \IZ$ by
\cite{Farrell-Jones(1995)}. 
We will not give the details of the proof of this $L$--theory statement.
The proof is however  completely analogous to the $K$--theory case. The extra 
assumption is needed to obtain a suitable compatibility with infinite products,
see~\ref{Lremark}.  The $L$--theory statements needed for this transition 
are provided in  \cite[Section 4]{Carlsson-Pedersen(1995)}.

I want to thank Tom Farrell, Wolfgang L{\"u}ck and Erik Pedersen for
helpful comments.


\section{Equivariant homology theories} \label{sec:orbit}

First let us briefly fix conventions on spectra. 
A spectrum $\IE$ is given by a sequence $(E_n)_{n \in \IN}$ of pointed spaces
and structure maps $\Sigma E_n \to E_{n+1}$. A map of spectra is a sequence of
maps $E_n \to F_n$ (for $n \in \IN$) that commutes with the structure maps.
A map of spectra is said to be a weak equivalence if it induces an 
isomorphism of (stable) homotopy groups.
Two spectra $\IE$ and $\IF$ are said to be weakly equivalent if there is
a zig-zag of weak equivalence
\[
\xymatrix
{
 \IE \ar[r]^\simeq & 
 \IA &
 \dots \ar[l]_\simeq \ar[r]^\simeq &
 \IF
}
\]
connecting $\IE$ to $\IF$.

Let $\Gamma$ be a group. 
The {\em Orbit Category} $\Or \Gamma$ has as objects the homogeneous 
spaces  $\Gamma/H$ and as morphisms $\Gamma$--equivariant maps 
$\Gamma/H \to \Gamma/K$ \cite{Bredon(1967)}. 
An $\Or \Gamma$--spectrum is a functor
from $\Or \Gamma$ to the category of spectra. A map of $\Or \Gamma$--spectra 
is a natural transformation. A map of $\Or \Gamma$--spectra is called 
a weak equivalence if it is a weak equivalence evaluated at every $\Gamma/H$.
Two $\Or \Gamma$--spectra are said to be weakly equivalent if they
are connected by a zig-zag of weak equivalences. 
Our main  example of an $\Or \Gamma$--spectrum is given by 
algebraic $K$--theory: for a ring $R$ there is an $\Or \Gamma$--spectrum
$\BFK R^{-\infty}$ whose value on $\Gamma/H$ is the $K$--theory spectrum of
the group ring $RH$. This functor has been constructed in 
\cite[Section 2]{Davis-Lueck(1998)}. In this paper we will denote spectra by
blackboard bold letters (like $\IE$) and $\Or \Gamma$--spectra by boldface
letters (like $\BFE$).

Associated to an $\Or \Gamma$--spectrum $\BFE$ is a
functor  from $\Gamma$--CW--complexes to spectra. 
Its value on a $\Gamma$--space 
$X$ is given by the {\em balanced smash product}  
\begin{eqnarray} \label{smashoveror}
 \IH^{\Or \Gamma}(X;\BFE) & 
 = & 
 X^H_+ \wedge_{\Or \Gamma} \BFE(\Gamma/H) \nonumber
 \\
 & = & 
 \coprod_{\Gamma/H} X^H_+ \wedge \BFE(\Gamma/H) / \sim, 
\end{eqnarray}
where $\sim$ is the equivalence relation generated by 
$(x\phi,y) \sim (x,\phi y)$
for $x \in X^K_+,y \in \BFE(\Gamma/H)$ and $\phi \co \Gamma/H \to \Gamma/K$
(cf.~\cite[Section 5]{Davis-Lueck(1998)}).
The homotopy groups of $\IH^{\Or \Gamma}(X;\BFE)$ will be denoted by
$H^{\Or \Gamma}_*(X;\BFE)$ and give an equivariant homology theory 
\cite[4.2]{Davis-Lueck(1998)}.   

A family of subgroups of $\Gamma$ is a collection of
subgroups of $\Gamma$ that is closed under conjugation and taking subgroups.
For such a family $\calf$ there is a classifying space $E\Gamma(\calf)$, 
namely a $\Gamma$--CW--complex characterized 
(up to $\Gamma$--homotopy equivalence) by the property 
that $E\Gamma(\calf)^H$ is contractible if 
$H \in \calf$ and empty otherwise. Given an $\Or \Gamma$--spectrum $\BFE$ there
is for any such family of subgroups $\calf$ the assembly map
$\IH^{\Or \Gamma}(E\Gamma(\calf);\BFE) \to 
 \IH^{\Or \Gamma}(\point;\BFE) = \BFE(\Gamma/\Gamma)$, 
cf.~\cite[Section 5]{Davis-Lueck(1998)}. This construction is natural in
the family $\calf$ and in this paper we will compare different families.

We will need the following recognition principle, 
cf.~\cite[6.3 2.]{Davis-Lueck(1998)}. A $\Gamma$--$\calf$--CW--complex,
is a $\Gamma$--CW--complex with isotropy groups contained in $\calf$.

\begin{lemma} \label{isoonorbits}
Let $\BFE \to \BFF$ be a map of $\Or \Gamma$--spectra. 
Let $\calf$ be a family of subgroups of $\Gamma$ such that
$\BFE(\Gamma/F) \to \BFF(\Gamma/F)$ is a weak equivalence for all
$F \in \calf$. Then  
\[
\IH^{\Or \Gamma} ( X;\BFE) \to \IH^{\Or \Gamma} (X ;\BFF)
\]
is a weak equivalence for any $\Gamma$--$\calf$--CW--complex.\qed
\end{lemma}

It will be useful for us to iterate the construction of $\Or \Gamma$--spectra,
i.e.\ define an $\Or \Gamma$--spectrum using the homology with respect to a 
different $\Or \Gamma$--spectrum.

\begin{lemma} \label{iterate}
Let $X,Y$ be $\Gamma$--CW--complexes and $\BFK$ be an $\Or \Gamma$--spectrum.
Define an $\Or \Gamma$--spectrum $\BFE$ by
\[
\BFE(\Gamma/H) = \IH^{\Or \Gamma}(\Gamma /H \x Y;\BFK).
\]
Then 
\[
\IH^{\Or \Gamma}(X;\BFE) \cong \IH^{\Or \Gamma}(X \x Y;\BFK).
\]
\end{lemma}

\proof
In the following formula $\Gamma/H$ will always correspond to the
first $\wedge_{\Or \Gamma}$ and $\Gamma/K$ to the second.
\begin{eqnarray*}
 \IH^{\Or \Gamma}(X;\BFE) & 
 = & 
 X^H_+ \wedge_{\Or \Gamma} \left( (\Gamma/H \x Y)^K_+ \wedge_{\Or \Gamma}
                       \BFK(\Gamma/K) \right)
 \\
 &
 = &
 \left( X_+^H \wedge_{\Or \Gamma} (\Gamma/H \x Y)^K_+ \right) 
                  \wedge_{\Or \Gamma} \BFK(\Gamma/K)
 \\
 &
 = &
 \left( \left(X^H_+ \wedge_{\Or \Gamma} (\Gamma/H)^K_+ \right) 
                  \wedge Y_+^K \right) 
                  \wedge_{\Or \Gamma} \BFK(\Gamma/K) 
 \\
 &
 \cong &
 (X^K_+ \wedge Y^K_+) \wedge_{\Or \Gamma} \BFK(\Gamma/K)
 \\
 &
 = &
 \IH^{\Or \Gamma}(X \x Y;\BFK).
\end{eqnarray*}
In the second, third and fourth line the first $\wedge_{\Or \Gamma}$ is a
balanced smash product with a  space, that is similarly defined 
as (\ref{smashoveror}). The homeomorphism from the third to  
the fourth line comes about as follows.
There is a natural $G$--action on $X^H_+ \wedge_{\Or \Gamma} (\Gamma / H)_+$
(where $G$ acts by multiplication on $\Gamma / H$, 
see \cite[7.1]{Davis-Lueck(1998)})
and by \cite[7.4.1]{Davis-Lueck(1998)} a natural $G$--homeomorphism
\[
X^H_+ \wedge_{\Or \Gamma} (\Gamma / H)_+ \cong X_+.
\]
Moreover, it is not hard to check that,
\[
X^H_+ \wedge_{\Or \Gamma} (\Gamma / H)_+^K = 
      (X^H_+ \wedge_{\Or \Gamma} (\Gamma / H)_+)^K.
\]  
Therefore,
$$
X^H_+ \wedge_{\Or \Gamma} (\Gamma / H)_+^K \cong X_+^K.
\eqno{\qed}$$

We finish this section with a formal splitting criterion.

\begin{proposition} \label{splitcrit}
Let $\BFE \to \BFF \to \BFG$ be maps of $\Or \Gamma$--spectra.
Let $\calf \subset \calg$ be families of subgroups of $\Gamma$.
Assume that $\BFE$ is weakly equivalent to 
$\Gamma/H \mapsto \IH^{\Or \Gamma}(\Gamma/H \x E\Gamma(\calf);\BFK)$
for some $\Or \Gamma$--spectrum $\BFK$. Assume moreover that 
$\BFE(\Gamma/F) \to \BFF(\Gamma/F)$ and
$\BFE(\Gamma/G) \to \BFG(\Gamma/G)$ 
are weak equivalences for all $F \in \calf$ and
$G \in \calg$.
Then
\[
H_*^{\Or \Gamma}(E\Gamma(\calf);\BFF) \to 
H_*^{\Or \Gamma}(E\Gamma(\calg);\BFF)
\]
is split injective.
\end{proposition}

\begin{proof}
Consider the following commutative diagram.
\[
\xymatrix
{
 \IH^{\Or \Gamma}(E\Gamma(\calf);\BFE) \ar[r]^{\alpha} \ar[d]^{\beta_0}  &
 \IH^{\Or \Gamma}(E\Gamma(\calg);\BFE) \ar[d] \arrow@/^5pc/[dd]^{\beta_1} 
 \\
 \IH^{\Or \Gamma}(E\Gamma(\calf);\BFF) \ar[r] \ar[d] &
 \IH^{\Or \Gamma}(E\Gamma(\calg);\BFF) \ar[d] 
 \\
 \IH^{\Or \Gamma}(E\Gamma(\calf);\BFG) \ar[r] &
 \IH^{\Or \Gamma}(E\Gamma(\calg);\BFG)  
}
\]
By the  first assumption and \ref{iterate} we have
\begin{eqnarray*}
\IH^{\Or \Gamma}(E\Gamma(\calf);\BFE) &
\simeq &
\IH^{\Or \Gamma}(E\Gamma(\calf) \x E\Gamma(\calf);\BFK),
\\
\IH^{\Or \Gamma}(E\Gamma(\calg);\BFE) & 
\simeq &
\IH^{\Or \Gamma}(E\Gamma(\calg) \x E\Gamma(\calf);\BFK).
\end{eqnarray*}
Now $\calf \subset \calg$ implies that both 
$E\Gamma(\calf) \x E\Gamma(\calf)$ and 
$E\Gamma(\calg) \x E\Gamma(\calf)$ are $\Gamma$--homotopy equivalent to
$E\Gamma(\calf)$. Thus $\alpha$ is a weak equivalence. The second
assumption and \ref{isoonorbits} imply that the maps labeled 
$\beta_i$ are also weak
equivalences.
\end{proof}


\section{Homotopy fixed points}
\label{sec:splitting}

A useful tool in proving injectivity results for assembly maps are homotopy
fixed points, cf.~\cite{Carlsson-Pedersen(1995)}.
Given an action of a group $\Gamma$ on a 
space $X$ the homotopy fixed points with respect to 
$\calf$ are by definition,
\[  
X^{h_\calf \Gamma} = \Map_{\Gamma}(E\Gamma(\calf),X).
\]
We will also need actions of $\Gamma$ on spectra. By definition $\Gamma$ acts
on a spectrum $\IE$, by acting (pointed) on each $E_n$ compatible with the 
structure maps. This allows to take (homotopy) fixed points level wise.
We will call a map $X \to Y$ a weak $\Or \Gamma$--equivalence,
if it is $\Gamma$--equivariant and induces a weak 
equivalence on all fixed point sets.

\begin{proposition} \label{assum}
Let
$\BFA,\BFB$ be $\Or \Gamma$--spectra with a $\Gamma$--action
(i.e.\ functors from $\Or \Gamma$ to spectra with $\Gamma$--action)
and $\calf \subset \calg$ two families of subgroups of $\Gamma$.
Assume that there is a
$\Gamma$--equivariant map of $\Or \Gamma$--spectra $\BFA \to \BFB$
such that the following holds.
\begin{enumerate}[\rm(1)] 
\item \label{assumhomology}
      There is an $\Or \Gamma$--spectrum $\BFK$ such that 
      the $\Or \Gamma$--spectra 
      $\BFA^\Gamma$ and 
      $\Gamma/H \mapsto \IH^{\Or \Gamma}(\Gamma/H \x E\Gamma(\calf);\BFK)$ 
      are weakly equivalent.
\item \label{assummap}
      For all  $G \in \calg$ there are 
      weak 
      $\Or \Gamma$--equivalences 
      \begin{eqnarray*}
      \BFA(\Gamma/G) &
      \simeq &
      \Map_G(\Gamma,\IA_0(G)) 
      \\
      \BFB(\Gamma/G) &
      \simeq &
      \Map_G(\Gamma,\IB_0(G))
      \end{eqnarray*}
      for 
      spectra $\IA_0(G),\IB_0(G)$ with a $G$--action. Moreover,
      there is a $G$--map $\IA_0(G) \to \IB_0(G)$ compatible with
      the $\Gamma$--map $\BFA(\Gamma / G) \to \BFB(\Gamma / G)$. 
\item \label{assumvircyc}
      For all $G \in \calg$ the induced map 
      $\IA_0(G)^G \to \IB_0(G)^{h_{\calf} G}$ 
      is a weak homotopy equivalence. (Here $\calf$ is viewed as the
      obvious family of subgroups of $G$ it induces.) 
\end{enumerate}
Then the map 
$H^{\Or \Gamma}_*(E\Gamma(\calf);\BFB^\Gamma) \to 
 H^{\Or \Gamma}_*(E\Gamma(\calg);\BFB^\Gamma)$ 
is split injective.
\end{proposition}

In our application in Section~\ref{sec:coefficient} $\calf$ will be the
family of finite subgroups and $\calg$ will be the family of 
virtual cyclic subgroups.
In order to prove \ref{assum}, we need three  lemmata.   
They will be used to relate 
fixed points of $\BFB$ (and $\BFA$) to homotopy fixed points of $\BFB$.
The proof of the first lemma is  straightforward.

\begin{lemma} \label{mapmap}
Let $H$ be a subgroup of $\Gamma$, $X$ a $\Gamma$--space and 
$Y$ an $H$--space. Then there is a natural homeomorphism
$$
\Map_\Gamma(X,\Map_H(\Gamma,Y)) \cong \Map_H(X,Y).
\eqno{\qed}$$
\end{lemma}

\begin{lemma} \label{fixhomotopyfix}
Let $H$ be a subgroup of $\Gamma$ and  $Y$ be an $H$--space. 
Let $S = \Map_H(\Gamma,Y)$. Then
\[
Y^H \cong S^\Gamma \quad \mbox{and} \quad  
   Y^{h_{\calf} H} \simeq  S^{h_\calf \Gamma}.
\]
If moreover $H \in \calf$ then
\[
S^\Gamma \simeq S^{h_\calf \Gamma}.
\]
\end{lemma}

\begin{proof}
Using \ref{mapmap} we have
\begin{eqnarray*}
S^\Gamma & = & \Map_\Gamma(\point,S) \\
         & = & \Map_\Gamma(\point,\Map_H(\Gamma,Y)) \\
         & \cong & \Map_H (\point , Y) \\
         & = & Y^H. \\
S^{h_{\calf}\Gamma} & = & \Map_\Gamma(E\Gamma(\calf),S) \\
                    & = & \Map_\Gamma(E\Gamma(\calf),\Map_H(\Gamma,Y)) \\
                    & \cong & \Map_H(E\Gamma(\calf),Y) \\
                    & \simeq & \Map_H(EH(\calf),Y) \\
                    & = & Y^{h_\calf H}.
\end{eqnarray*}
To prove the last assertion, observe that
if $H \in \calf$, then $EH(\calf)$ is a point and $Y^{h_\calf H} = Y^H$.
Therefore $S^{h_\calf \Gamma} \simeq S^\Gamma$.
\end{proof}

\begin{lemma} \label{coefiso}
For $F \in \calf$ and $G \in \calg$ the induced maps 
\begin{eqnarray*}
\BFB(\Gamma/F)^\Gamma &
\to &
\BFB(\Gamma/F)^{h_\calf \Gamma},
\\
\BFA(\Gamma/G)^\Gamma &
\to &
\BFB(\Gamma/G)^{h_\calf \Gamma}
\end{eqnarray*}
are homotopy equivalences.
\end{lemma}

\begin{proof}
The first homotopy equivalence  follows easily 
from \ref{assum} (\ref{assummap}) and the 
second part of \ref{fixhomotopyfix}. The second map is by
\ref{assum} (\ref{assummap}) and the first part of \ref{fixhomotopyfix}  
equivalent to $\IA_0(G)^G \to \IB_0(G)^{h_\calf G}$ and a 
homotopy equivalence by \ref{assum} (\ref{assumvircyc}).
\end{proof}

\begin{proof}[Proof of Proposition~\ref{assum}]
Set $\BFE = \BFA^{\Gamma}$, $\BFF = \BFB^{\Gamma}$ and
$\BFG = \BFB^{h_\calf \Gamma}$. 
In order to apply \ref{splitcrit}, we need to check
that $\BFA(\Gamma/G)^{\Gamma} \to \BFB(\Gamma/G)^{h_\calf \Gamma}$
and $\BFA(\Gamma/F)^{\Gamma} \to \BFB(\Gamma/F)^{\Gamma}$ 
are weak equivalences for
$G \in \calg$ and $F \in \calf$. This a consequence of \ref{coefiso}.
\end{proof}


\section{Controlled algebra}

Let $Z$ be a topological space and $R$ be a ring.
Controlled algebra is concerned with categories of
$R$--modules over $Z$ 
($M=\bigoplus_{z \in Z} M_z$) and $R$--module maps over $Z$ 
($\phi = (\phi_{z,z'} \co M_z' \to M_z)$). 
We will need an equivariant version of this theory that has 
been studied in \cite{Bartels-Farrell-Jones-Reich(2001)}.
Let $\Gamma$ be a group and $X$ be a $\Gamma$--space. 
The equivariant continuous control condition 
$\cale_{\Gamma cc}(X)$ (consisting of subsets of 
$(X \x \einsu)^{\x 2}$) is defined in 
\cite[2.5]{Bartels-Farrell-Jones-Reich(2001)}. Let 
$p \co Y \to X$ be a continuous $\Gamma$--map. We define a 
category $\calc(Y;p)$ of $R$--modules over 
$Y \x \Gamma \x \einsu$: Its objects are locally finite 
(see \cite[Section 2.2]{Bartels-Farrell-Jones-Reich(2001)}) 
free $R$--modules $M = \bigoplus M_{(y,\gamma,t)}$ 
subject to the condition that there is a compact subset 
$K \subset Y \x \Gamma$ (depending on $M$) such that 
$M_{(y,\gamma,t)}=0$ unless $(y,\gamma) \in \Gamma K$. 
Morphisms $\phi=(\phi_{(y,\gamma,t),(y',\gamma',t')})$ 
are required to satisfy the following condition: 
there is $E \in \cale_{\Gamma cc}(X)$ (depending on $\phi$)
such that $\phi_{(y,\gamma,t),(y',\gamma',t')}=0$ unless 
$((p(y),t),(p(y'),t')) \in E$. Note that this definition 
depends on the group action we have in mind. 
The objects of the full subcategory 
$\calc_0(Y;p) \subset \calc(Y;p)$  
have by definition support in $Y \x \Gamma \x [1,\alpha]$, 
i.e.\ for every module $M$ 
there  is $\alpha > 0$ such that $M_{y,\gamma,t} = 0$ unless $t \leq \alpha$.
This inclusion is a Karoubi filtration 
(\cite[1.27]{Carlsson-Pedersen(1995)}) and
we denote the quotient by $\cald(Y;p)$. The group
$\Gamma$ acts on all these categories. The fixed point category 
$\cald^\Gamma(Y;p)$ appeared in 
\cite{Bartels-Farrell-Jones-Reich(2001)}. We abbreviate
\[
\IK(p) = \IK^{-\infty} \cald(Y;p).
\]
If $p = \id_X$ we will write $\IK(X)$ for $\IK(\id_X)$. 
An important application of controlled algebra has been the construction
of homology theories \cite{Pedersen-Weibel(1989)}.
The following equivariant version of this result is proven in 
\cite[Section 5 and 6.2]{Bartels-Farrell-Jones-Reich(2001)}.

\begin{theorem} \label{homology}
The functor 
\begin{eqnarray*} 
X \mapsto \Omega\IK(X)^\Gamma
\end{eqnarray*}
from $\Gamma$--CW--complexes to spectra
is weakly equivalent to
$$
X \mapsto \IH^{\Or \Gamma}(X;\BFK R^{-\infty}).
\eqno{\qed}$$
\end{theorem}

We will later on need the following simple observation.

\begin{lemma} \label{point}
$\IK(X \x Y \to Y) \to \IK(Y)$ is a weak $\Or \Gamma$--homotopy equivalence.
\end{lemma}

\begin{proof}
It is not hard to check that 
$\cald^H(X \x Y;X \x Y \to Y)) \to 
 \cald^H(Y;\id_Y)$ is an equivalence of categories
for any subgroup $H$.
\end{proof}

The next lemma will later on be the key ingredient in checking  
condition~\ref{assum} (\ref{assummap}).

\begin{lemma} \label{MapandGoverH}
Let $p \co X \to Y \x \Gamma/H$ be a $\Gamma$--map. Let 
$X_0 = p^{-1}(Y \x \{eH\} )$
and denote by $p^H_0 \co X_0 \to Y$ the $H$--map induced by $p$. Then there is
a weak $\Or \Gamma$--equivalence
\[
\IK(p) \simeq \Map_H(\Gamma,\IK(p^H_0)).
\]
\end{lemma}

\proof
For $U \subset \Gamma / H$
let $X [ U ] = p^{-1}(Y \x U)$ and 
$p [ U ] = p|_{X[U]}$.
For a subgroup $F$ we abbreviate 
$\calc^F [ U ] = \cald^F (X [ U ] ; p [ U ])$.
Clearly $\IK(p_0^H) = \IKi \calc [ e H ]$. The continuous control condition 
$\cale_{\Gamma cc}(Y \x \Gamma/H)$
separates in particular different path components. 
Therefore we get 
\[
\cald(X;p) \cong \prod_{\gamma H \in \Gamma / H} \calc [ \gamma H ].
\]
Projections induce a map 
\[
\IKi  \cald(X;p) \to 
 \prod_{\gamma H \in \Gamma / H} \IKi \calc [ \gamma H ] 
            \cong \Map_H(\Gamma,\IK(p^H_0)).
\]
We have to show that this map is a weak $\Or \Gamma$--equivalence. 
Let $F$ be a subgroup of $\Gamma$. Again, the continuous control condition
implies 
\[
\cald^F(X;p) \cong \prod_{F \gamma H \in F \backslash \Gamma /H} 
                   \calc^F [ F \gamma H ].
\]
Using the fact that $\IKi$ commutes with fixed points and up to weak
equivalence with infinite products 
\cite{Carlsson(1995)} we obtain
\begin{eqnarray*}
(\IKi  \cald(X;p))^F & \cong & \IKi \cald^F(X;p) \\ 
   &  \simeq  & \prod_{F \gamma H \in F \backslash \Gamma /H} 
                    \IKi \calc^F [ F \gamma H ].
\end{eqnarray*}
Moreover, 
\begin{eqnarray*}
\IKi \calc^F [ F \gamma H ] & 
     \cong & \IKi \calc^{F \cap \gamma H \gamma^{-1}} [ \gamma H ] \\ 
 & \cong & (\IKi  \calc [ \gamma H ] )^{F \cap \gamma H \gamma^{-1}} \\
 &  \cong  & \Big(\prod_{f \gamma H \in  F (\gamma H)} 
                                \IKi \calc [ f \gamma H ] \Big)^F. \\
\end{eqnarray*}
(Here $F(\gamma H)$ denotes the $F$--orbit of $\gamma H$ in $\Gamma / H$.)
We finish the argument by observing that
$$
\Big( \prod_{\gamma H \in \Gamma / H}  \IKi \calc [ \gamma H ] \Big)^F 
\cong 
\prod_{F \gamma H \in F \backslash \Gamma /H} 
  \Big( \prod_{f \gamma H \in F (\gamma H)} \IKi \calc [ f \gamma H ] 
                 \Big)^F.\eqno{\qed}$$

\begin{remark} \label{Lremark}
In the proof above we used the compatibility of $K$--theory 
with infinite products
from \cite{Carlsson(1995)}. At this point the $L$--theory version of
our splitting result needs the additional 
assumption stated in the introduction.
It is explained in \cite[p. 756]{Carlsson-Pedersen(1995)}
that for additive categories with involutions $\cala_n$ 
there is a weak equivalence  
\[ 
\ILi \Big( \prod \cala_n \Big) \cong \prod \ILi \cala_n,
\]
provided there is $i_0$ independent of $n$ such that 
$K_{-i} \cala_n = 0$ for all $i \geq i_0$. Thus, an $L$--theory 
version of \ref{MapandGoverH} needs an additional assumption. 
A sufficient assumption is that $K_{-i} RH = 0$ for all sufficiently large $i$.
\end{remark}

Under sufficient control conditions, 
there is no difference between fixed points
and homotopy fixed points. This is an important ingredient in the proof
of injectivity of assembly maps in \cite{Carlsson-Pedersen(1995)}
and \cite{Rosenthal(2002)}. We will 
need the following version of this result.

\begin{lemma} \label{descent}
Let $X$ be a cocompact  $\Gamma$--CW--complex with isotropy groups 
contained in a family 
of subgroups $\calf$. Then the obvious map
\[
\IK(X)^\Gamma \to \IK(X)^{h_\calf \Gamma}
\]
is a homotopy equivalence.
\end{lemma}

\begin{proof}
This is \cite[6.2]{Rosenthal(2002)}.
One proceeds by induction on the 
equivariant cells of $X$. The induction step uses  
\ref{MapandGoverH} and \ref{fixhomotopyfix}.
\end{proof}

The following result is closely related to  
\cite[7.1]{Rosenthal(2002)}. Using what is sometimes 
called the descent principle 
it can be used to show split injectivity of (\ref{ktheorysplit}) in the
base case, i.e.\ for virtually cyclic $\Gamma$. 
(The point of the 
descent principle is that it requires 
only knowledge about fixed points of finite subgroups.)
The infinite cyclic and
the infinite dihedral group act properly on $\IR$.
Virtually cyclic groups map either onto the integers or the infinite dihedral
group (\cite[2.5]{Farrell-Jones(1995)}), 
and act therefore also properly on $\IR$. 
The restriction of this action to
finite subgroups is either trivial or factors through the action of $\IZ/2$ by
a reflection.

\begin{proposition} \label{finassemblyiso}
Consider $\IR$ with the aforementioned proper 
action of a virtual cyclic group $V$. 
If $H$ is a finite subgroup of $V$, then 
\[
\IK(\IR)^H \to \IK(\IR \to \point)^H
\]
is a weak equivalence.
\end{proposition}

In order to prove this we will need a slightly different 
construction of $\cald(Y;p)$ for a continuous $\Gamma$--map
$p \co Y \to X$ where $X$ carries a $\Gamma$--equivariant metric $d$.
Define
the subcategory $\tilde\calc(Y;p) \subset \calc(Y;p)$ whose morphisms have to
satisfy the additional condition, that there is $\alpha > 0$
(depending on $\phi$) such that 
$\phi_{(y,\gamma,t),(y',\gamma',t')}=0$  unless 
$d(p(y),p(y')) \leq \alpha$. The corresponding inclusion
$\tilde\calc_0(Y;p) \subset \tilde\calc(Y;p)$ is again a Karoubi filtration.
It is not to hard to check, that its quotient $\tilde\cald(Y;p)$
is equivalent to  $\cald(Y;p)$ and that this is  compatible with the
$\Gamma$--actions, cf.~\cite[8.8]{Bartels-Farrell-Jones-Reich(2001)}. 
However, one has to be a little careful with the definitions
to get this even before taking fixed points. 
In particular, it is at this point  important that
all $E \in \cale_{\Gamma cc}(X)$ are required to be $\Gamma$-invariant,
\cite[2.5(iii)]{Bartels-Farrell-Jones-Reich(2001)}.

\begin{lemma} \label{squeezing}          
The $K$--theory of $\tilde\calc^H(\IR;\id_\IR)$ vanishes
under the assumption of \ref{finassemblyiso}. 
(Here we consider the standard metric on $\IR$.)
\end{lemma} 

\begin{proof}
Let $x_0 \in \IR$ be a fixed point for the action of $H$.
We will need various full subcategories of $\tilde\calc^H(\IR;\id_\IR)$.
Let $\tilde\cals$ 
be the full subcategory whose objects 
have support in $[x_0 - \alpha,x_0 + \alpha] \x V \x \einsu$ 
            for some $\alpha > 0$;
$\tilde\cals_{+}$ 
be the full subcategory whose objects 
have support in $[x_0,x_0+\alpha] \x V \x \einsu$ for some $\alpha > 0$;
$\tilde\cals_{-}$ 
be the full subcategory whose objects 
have support in $[x_0-\alpha,x_0] \x V \x \einsu$ for some $\alpha > 0$;
$\tilde\calc_{+}$
be the full subcategory whose objects 
have support in $[x_0,\infty) \x V \x \einsu$;
$\tilde\calc_{-}$
be the full subcategory whose objects 
have support in $(-\infty,x_0] \x V \x \einsu$.
Then $\tilde\cals \subset \tilde\calc^H(\IR;\id_\IR)$,
$\tilde\cals_{+} \subset \tilde\calc_{+}$ and
$\tilde\cals_{-} \subset \tilde\calc_{-}$
are Karoubi filtrations and we denote the quotient categories by
$\tilde\calq$, $\tilde\calq_+$ and $\tilde\calq_-$. It is not hard to 
check that the first of these quotients is equivalent to the direct sum
of the two later. The $K$--theory of $\tilde\calq$ is therefore the sum of
the $K$--theories of $\tilde\calq_+$ and $\tilde\calq_-$.
Applying $\IKi$ to Karoubi filtrations 
gives a homotopy fibration by \cite[1.28]{Carlsson-Pedersen(1995)}.
Putting all this together, we see that it 
suffices to show that the $K$--theory of each of 
our five full subcategories is trivial. 
The map $(x,v,t) \mapsto ((x-x_0)/2 + x_0,v,t+1)$ induces an 
Eilenberg swindle on $\tilde{\cals}$, $\tilde{\cals_+}$ and $\tilde{\cals_-}$;
the maps $(x,v,t) \mapsto (x+1,v,t)$ and $(x,v,t) \mapsto (x-1,v,t)$
induce Eilenberg swindles on $\tilde{\calc_+}$ and $\tilde{\calc_-}$.
\end{proof}

Note that it is important to use the   
category $\tilde\calc$ rather than $\calc$ for this argument.
For example, the corresponding subcategory $\cals$ of 
$\calc^H(\IR;\id_\IR)$ is not a Karoubi filtration.

\begin{proof}[Proof of \ref{finassemblyiso}]
Let $p$ denote the projection $\IR \to \point$. We will use the following 
diagram.
\[
\xymatrix
{
\tilde\calc_0(\IR;\id_\IR)^H \ar[r] \ar[d]_{\calf_1} &
\tilde\calc(\IR;\id_\IR)^H \ar[r] \ar[d]_{\calf_2} &
\tilde\cald(\IR;\id_\IR)^H  \ar[d]_{\calf_3} 
\\
\tilde\calc_0(\IR;p)^H \ar[r]  &
\tilde\calc(\IR;p)^H \ar[r]  &
\tilde\cald(\IR;p)^H  
\\
}
\]
It is not hard to check that $\calf_1$ is an equivalence of categories. 
The $K$--theory of $\tilde\calc(\IR;\id_\IR)^H$ vanishes by \ref{squeezing}.
The map $(x,v,t) \mapsto (x,v,t+1)$ gives an Eilenberg swindle on 
$\tilde\calc(\IR;p)^H$
and its $K$--theory also vanishes. As used before,
applying $\IKi$ to Karoubi filtrations 
gives a homotopy fibration by \cite[1.28]{Carlsson-Pedersen(1995)}. Thus 
$\calf_3$ induces an isomorphism in $K$--theory. 
The result follows, since 
$\tilde\cald(\IR;q) = \cald(\IR;q)$ 
for any $q$ as noted before \ref{squeezing}. 
\end{proof}


\section{The coefficient spectra}
\label{sec:coefficient}

This section contains the proof of Theorem~\ref{main} from the 
introduction.
As before, we fix a ring $R$ and a group $\Gamma$.
For a subgroup $H$ of $\Gamma$ let 
\[
p_{\Gamma/H} \co 
\Gamma/H \x E\Gamma(\calfin) 
\to \Gamma/H
\]
be the obvious projections. 
We define two $\Or \Gamma$ spectra $\BFA$ and $\BFB$ by
\begin{eqnarray*}
\BFA(\Gamma/H) &
= &
\IK(\Gamma /H \x E\Gamma(\calfin)),
\\
\BFB(\Gamma/H) &
= &
\IK(p_{\Gamma/H}).
\end{eqnarray*}
Both, $\BFA$ and $\BFB$ are naturally equipped with a $\Gamma$--action.
There is an obvious $\Gamma$--equivariant map of $\Or\Gamma$--spectra 
$\BFA \to \BFB$.

We will show 
that these spectra satisfy the hypothesis of \ref{assum}
with respect to the families $\calfin \subset \calvc$.
For \ref{assum} (\ref{assumhomology}) this follows from  
\ref{homology}, where $\BFK$ is the 
algebraic $K$--theory $\Or \Gamma$--spectrum
$\BFK R^{-\infty}$. In \ref{coind} we will prove 
that \ref{assum} (\ref{assummap}) is satisfied. 
The final condition
\ref{assum} (\ref{assumvircyc}) will follow from \ref{einsokay}.
Moreover, it is an easy consequence of \ref{homology} and \ref{point}
that $\Omega \BFB^\Gamma$ is weakly equivalent to $\BFK R^{-\infty}$
and therefore Theorem~\ref{main} will be a consequence of the splitting 
result \ref{assum}.
  
For a subgroup $H$ of $\Gamma$ let 
\begin{eqnarray*}
\IA_0(H) &
= &
\IK(\res_\Gamma^H E\Gamma(\calfin)),
\\
\IB_0(H) &
= &
\IK(\res_\Gamma^H(E\Gamma(\calfin) \to \point)).
\end{eqnarray*}
Here $\res_\Gamma^H$ denotes the forgetful functor from $\Gamma$--spaces
to $H$--spaces.

The next statement is an immediate consequence of \ref{MapandGoverH}
and verifies \ref{assum} (\ref{assummap}).

\begin{lemma} \label{coind}
There are natural weak $\Or \Gamma$ equivalences
\begin{align*}
\BFA(\Gamma/H)\ &
=\ 
\Map_H(\Gamma,\IE_0(H)),
\\
\BFB(\Gamma/H)\ &
=\ 
\Map_H(\Gamma,\IB_0(H)).\tag*{\qed}
\end{align*}
\end{lemma}

Finally, we verify \ref{assum} (\ref{assumvircyc}).

\begin{proposition} \label{einsokay}
For $V \in \calvc$ the obvious map  
\[
\xymatrix
{
\IA_0(V)^V = \IK(\res_\Gamma^V E\Gamma(\calfin))^V  \ar[d] 
\\
\IB_0(V)^{h_\calf V} = 
    \IK(\res_\Gamma^H(E\Gamma(\calfin) \to \point))^{h_\calfin V}  
}
\]
is a homotopy equivalence.
\end{proposition}

\begin{proof} 
We can choose $EV(\calfin) = \IR$ with the proper action
used towards the end of the previous section.
We will use the following commutative diagram.
\[
\xymatrix
{
\IK(\IR)^V \ar[rr]^{\alpha_0} \ar[d]^{\alpha_1} & &
\IA_0(V)^V  \ar[dd]
\\
\IK(\IR)^{h_\calfin V} \ar[d]^{\beta_0}
\\
\IK(\IR \to \point)^{h_\calfin V} \ar[rr]^{\beta_1} & &
\IB_0(V)^{h_\calfin V} &
}
\]
The maps labeled $\alpha_i$ and $\beta_i$ are all homotopy equivalences:
$\alpha_0$ by the fact that $\res_\Gamma^V E\Gamma(\calfin))$ is 
also an $EV(\calfin)$ and \ref{homology} and $\alpha_1$ by \ref{descent}.
To study the maps labeled $\beta_i$ we need a fact about homotopy fixed points:
if an equivariant map induces a homotopy equivalence 
on fixed points for finite subgroups,
then it induces a homotopy equivalence on homotopy fixed points with 
respect to $\calfin$,
see \cite[4.1]{Rosenthal(2002)}. 
Thus $\beta_1$ is a homotopy equivalence by \ref{point}.
The map $\IK(\IR) \to \IK(\IR \to \point)$ 
induces a homotopy equivalence on fixed points under all finite 
subgroups of $V$ by \ref{finassemblyiso} and therefore $\beta_0$ 
is also a homotopy equivalence.
\end{proof}


{\makeatletter
\@thebibliography@{Bre67}
\makeatother

\bibitem[BFJR]{Bartels-Farrell-Jones-Reich(2001)}
A.~Bartels, F.~T. Farrell, L.~E. Jones, and H.~Reich.
\newblock On the isomorphism conjecture in algebraic ${K}$--theory.
\newblock to appear in {\em Topology}.

\bibitem[BHM93]{Boekstedt-Hsiang-Madsen(1993)}
M.~B{\"o}kstedt, W.~C. Hsiang, and I.~Madsen.
\newblock The cyclotomic trace and algebraic ${K}$--theory of spaces.
\newblock {\em Invent. Math.}, 111(3):465--539, 1993.

\bibitem[Bre67]{Bredon(1967)}
G.~E. Bredon.
\newblock {\em Equivariant cohomology theories}.
\newblock Springer-Verlag, Berlin, 1967.

\bibitem[Car95]{Carlsson(1995)}
G.~Carlsson.
\newblock On the algebraic ${K}$--theory of infinite product categories.
\newblock {\em $K$--Theory}, 9(4):305--322, 1995.

\bibitem[CP95]{Carlsson-Pedersen(1995)}
G.~Carlsson and E.~K. Pedersen.
\newblock Controlled algebra and the {N}ovikov conjectures for ${K}$-- and
  ${L}$--theory.
\newblock {\em Topology}, 34(3):731--758, 1995.

\bibitem[DL98]{Davis-Lueck(1998)}
J.~F. Davis and W.~L{\"u}ck.
\newblock Spaces over a category and assembly maps in isomorphism conjectures
  in ${K}$-- and ${L}$--theory.
\newblock {\em $K$--Theory}, 15(3):201--252, 1998.

\bibitem[FJ93]{Farrell-Jones(1993)}
F.~T. Farrell and L.~E. Jones.
\newblock Isomorphism conjectures in algebraic ${K}$--theory.
\newblock {\em J. Amer. Math. Soc.}, 6(2):249--297, 1993.

\bibitem[FJ95]{Farrell-Jones(1995)}
F.~T. Farrell and L.~E. Jones.
\newblock The lower algebraic {$K$}--theory of virtually 
           infinite cyclic groups.
\newblock {\em $K$--Theory}, 9(1):13--30, 1995.

\bibitem[FJ98]{Farrell-Jones(1998)}
F.~T. Farrell and L.~E. Jones.
\newblock Rigidity for aspherical manifolds with 
     $\pi\sb 1\subset {G}{L}\sb m( {R} )$.
\newblock {\em Asian J. Math.}, 2(2):215--262, 1998.

\bibitem[L{\"u}02]{Lueck(2000)}
W.~L{\"u}ck.
\newblock Chern characters for proper equivariant homology theories and
  applications to ${K}$-- and ${L}$--theory.
\newblock {J. Reine Angew. Math.}, 543:193--234, 2002.

\bibitem[PW89]{Pedersen-Weibel(1989)}
E.~K. Pedersen and C.~A. Weibel.
\newblock ${K}$--theory homology of spaces.
\newblock In {\em Algebraic topology (Arcata, CA, 1986)}, pages 346--361.
  Springer, Berlin, 1989.

\bibitem[Ros03]{Rosenthal(2002)}
D.~Rosenthal.
\newblock {\em Splitting with continuous control in algebraic $K$--theory}.
\newblock {\tt arXiv:math.AT/0309106}

\endthebibliography}

\Addresses\recd

\end{document}

%% file: agtout.tex
\def\ifplaintex{\expandafter\ifx\csname documentclass\endcsname\relax}

\def\gtp{{\mathsurround=0pt\it $\cal G\mskip-2mu$eometry \&\ 
$\cal T\!\!$opology $\cal P\!$ublications}}  

\def\recd{{\small Received:\qua\receiveddate\ifx\reviseddate\relax
\else\qquad Revised:\qua\reviseddate\fi\par}} 

\def\lognumber#1{\def\thelognumber{#1}}
\def\volumenumber#1{\def\thevolumenumber{#1}}
\def\volumeyear#1{\def\thevolumeyear{#1}}
\def\papernumber#1{\def\thepapernumber{#1}}
\def\pagenumbers#1#2{\def\startpage{#1}\def\finishpage{#2}}
\def\published#1{\def\publishdate{#1}}

\def\received#1{\def\receiveddate{#1}}
\def\revised#1{\def\reviseddate{#1}}
\def\accepted#1{\def\accepteddate{#1}}
\def\asciititle#1{\def\theasciititle{#1}}

\def\asciiaddress#1{\def\theasciiaddress{#1}}

\long\def\asciiabstract#1{\long\def\theasciiabstract{#1}}
\def\asciikeywords#1{\def\theasciikeywords{#1}}

\let\\\par\let\thelognumber\relax\let\thevolumenumber\relax
\let\thepapernumber\relax\let\thevolumeyear\relax\let\startpage\relax
\let\finishpage\relax\let\publishdate\relax\let\receiveddate\relax
\let\reviseddate\relax\let\accepteddate\relax\let\theasciititle\relax
\let\theasciiauthors\relax\let\theasciiaddress\relax
\let\theasciiabstract\relax\let\theasciikeywords\relax

\let\theasciiemail\relax

\ifplaintex
\font\logobig=cmssbx10 scaled 3836
\font\logomed=cmssbx10 scaled 2557
\else
\font\logobig=cmssbx10 scaled 4200
\font\logomed=cmssbx10 scaled 2800
\fi

\long\def\makeagttitle{   
\count0=\startpage
\agt\hfill      
\hbox to 45truept{\vbox to 0pt{\vglue -13truept{\logomed A\kern -.37em{\logobig 
T}\kern -.38em G}\vss}\hss}
\break
{\small Volume \thevolumenumber\ (\thevolumeyear)
\startpage--\finishpage\nl
Published: \publishdate}

\vglue .25truein

{\parskip=0pt\leftskip 0pt plus
1fil\def\\{\par\smallskip}{\Large\bf\thetitle}\par\medskip} \vglue
0.05truein

{\parskip=0pt\leftskip 0pt plus 1fil\def\\{\par}{\sc\theauthors}
\par\medskip}%
 
\vglue 0.03truein 

{\small\leftskip 25truept\rightskip 25truept{\bf Abstract}\stdspace\theabstract

{\bf AMS Classification}\stdspace\theprimaryclass
\ifx\thesecondaryclass\relax\else; \thesecondaryclass\fi\par
{\bf Keywords}\stdspace \thekeywords\par}\vglue 7truept

}   

\ifplaintex
\font\phead=cmsl9 scaled 950
\font\pnum=cmbx10 scaled 913
\font\pfoot=cmsl9 scaled 950
\headline{\vbox to 0pt{\vskip -4.5mm\line{\small\phead\ifnum
\count0=\startpage ISSN 1472-2739 (on-line) 1472-2747 (printed)
\hfill {\pnum\folio}\else\ifodd\count0\def\\{ }%
\ifx\theshorttitle\relax\thetitle\else\theshorttitle\fi\hfill{\pnum\folio}
\else\def\\{ and }{\pnum\folio}\hfill\ifx\theshortauthors\relax\theauthors
\else\theshortauthors\fi\fi\fi}\vss}}
\footline{\vbox to 0pt{\vglue 0mm\line{\small\pfoot\ifnum\count0=\startpage
\copyright\ \gtp\hfill\else
\agt, Volume \thevolumenumber\ (\thevolumeyear)\hfill\fi}\vss}}
\else
\font=cmsl9 scaled 1050
\font\lnum=cmbx10 
\font=cmsl9 scaled 1050
\makeatletter
\def\@oddhead{{\small\lhead\ifnum\count0=\startpage ISSN 1472-2739 
(on-line) 1472-2747 (printed)\hfill {\lnum\number\count0}\else\ifodd\count0
\def\\{ }\ifx\theshorttitle\relax \thetitle \else\theshorttitle\fi\hfill
{\lnum\number\count0}\else\def\\{ and }{\lnum\number\count0}
\hfill\ifx\theshortauthors\relax 
\theauthors\else\theshortauthors\fi\fi\fi}}\def\@evenhead{\@oddhead}
\def\@oddfoot{\small\lfoot\ifnum\count0=\startpage\copyright\ \gtp\hfill\else
\agt, Volume \thevolumenumber\ (\thevolumeyear)\hfill\fi}
\def\@evenfoot{\@oddfoot}
\makeatother
\fi
\let\maketitlepage\makeagttitle
\let\makeshorttitle\maketitlepage
\let\maketitle\maketitlepage

\newwrite\gtoutfile
\long\gdef\makeheadfile{  
{\def\\{, }\def\s{ }
\immediate\openout\gtoutfile head.xxx
\immediate\write\gtoutfile{Proxy-for: \ifx\theasciiauthors\relax
\theauthors\else\theasciiauthors\fi\s<\ifx\theasciiemail\relax\theemail\else\theasciiemail\fi>}
\immediate\write\gtoutfile{\noexpand\\}
\immediate\write\gtoutfile{Authors: \ifx\theasciiauthors\relax
\theauthors\else\theasciiauthors\fi}
{\def\\{ }\immediate\write\gtoutfile{Title: \ifx\theasciititle\relax
\thetitle\else\theasciititle\fi}}
\immediate\write\gtoutfile{Subj-class: GT or SG, GR etc}
\immediate\write\gtoutfile{MSC-class: \theprimaryclass\ifx\thesecondaryclass\relax\else, \thesecondaryclass\fi}
\immediate\write\gtoutfile{Journal-ref: Algebr. Geom. Topol. \thevolumenumber\s
(\thevolumeyear) \startpage-\finishpage}
\immediate\write\gtoutfile{Comments: Published by Algebraic and
Geometric Topology at}
\immediate\write\gtoutfile{\s\s\s  http://www.maths.warwick.ac.uk/agt/AGTVol\thevolumenumber/agt-\thevolumenumber-\thepapernumber.abs.html}
\immediate\write\gtoutfile{\noexpand\\}
\immediate\write\gtoutfile{}
\ifx\theasciiabstract\relax
\immediate\write\gtoutfile{\theabstract}\else
\immediate\write\gtoutfile{\theasciiabstract}\fi
\immediate\write\gtoutfile{}
\immediate\write\gtoutfile{\noexpand\\}
\immediate\write\gtoutfile{}
\immediate\closeout\gtoutfile}}  

\def\maketitlepage{\makeagttitle\makeheadfile}
\let\makeshorttitle\maketitlepage
\let\maketitle\maketitlepage